\documentclass[11pt]{article}
\usepackage[body={7.0in,9.0in}]{geometry}
\usepackage{systeme}
\usepackage[dvips]{color}
\usepackage{amssymb}
\usepackage{amsmath}
\usepackage{amsthm}
\usepackage{hyperref}
\usepackage{algorithm2e}

\def\F{\mathbb F}

\def\S{\mathbb S}

\def\A{\mathbb A}
\def\N{\mathbb N}
\def\O{\mathbb O}

\ifx\pdfoutput\undefined
\usepackage{graphicx}
\else
\usepackage[pdftex]{graphicx}
\fi

\def\d{\underline{d}}

\def\0{\underline{0}}

\usepackage{tikz}
\usepackage{amsfonts}

\newtheorem{theorem}{Theorem}[section]
\newtheorem{lemma}{Lemma}[section]
\newtheorem{corollary}{Corollary}[section]
\newtheorem{proposition}{Proposition}[section]
\newtheorem{example}{Example}[section]
\newtheorem{remark}{Remark}[section]
\newtheorem{definition}{Definition}[section]

\begin{document}
\begin{center}
{\fontsize{14}{20}\bf
Enumeration of irreducible and extended irreducible Goppa codes
}
\end{center}

\begin{center}
\textbf{Kondwani Magamba$^{1,2}$ and John A. Ryan$^3$}\\[.25cm]
$^1$Malawi University of Science and Technology, Malawi\\ 
$^{2}$Mzuzu University, Malawi\\
$^{3}$Chombe Boole Research Center, Malawi

\end{center}

\begin{abstract}
We obtain upper bounds on the number of irreducible and extended irreducible Goppa codes over $\F_p$ of length $q$ and $q+1$, respectively defined by polynomials of degree $r$, where $q=p^t$ and $r\geq 3$ is a positive integer.
\end{abstract} 

\section{Introduction}
It is well known that Goppa codes have few invariants and that the number of inequivalent Goppa codes grows exponentially with the length and dimension of the code \cite{Ryan_Disc}. These facts led to the exploitation of Goppa codes in the McEliece cryptosystem. There has been research work on the enumeration of extended Goppa codes but most of it has been confined to particular cases. The reference \cite{Ryan4} gives an upper bound on the number of inequivalent extended irreducible binary quartic Goppa codes. Recently, an upper bound on the number of inequivalent extended irreducible Goppa codes of degree $r$ and length $q^n+1$, where $q=p^t$ with the restriction that $n$ and $r$ be primes was found in \cite{Khowa2}. Also, in 1978, Chen \cite{Chen} gave upper bounds on the number of inequivalent irreducible and extended irreducible Goppa codes of length $q^n+1$ which are not tight. In this paper we derive upper bounds on the number of irreducible and extended irreducible Goppa codes which are tighter than the bounds found in \cite{Chen}. Our approach takes advantage of recent work by various researchers on the action of $PGL(2,q)$ and $GL(2,q)$ on the set of irreducible polynomials in $\F_q[x]$, see \cite{Reis} and \cite{Gare}. This work sheds more light on the structure of Goppa codes and the strength of the McEliece cryptosytem.

\section{Preliminaries}
\subsection{Irreducible and extended irreducible Goppa Codes}
We begin by defining an irreducible Goppa code.
\begin{definition}
Let $n$ be a positive integer, $q$ be a power of a prime number $p$ and $g(z)\in \F_{q^n}[z]$ be irreducible of degree $r$. Let
$L=\F_{q^n}=\{\zeta_i: 0\leq i \leq q^n-1 \}$. Then an irreducible Goppa code $\Gamma(L,g)$ is defined as the set of all vectors $\underline{c}=(c_0,c_1,\ldots,c_{q^n-1})$ with components in $\F_q$ which satisfy the condition
\begin{equation}
\sum_{i=0}^{q^n-1}\frac{c_i}{z-\zeta_i}\equiv 0~ \mbox{mod}~g(z). 
\end{equation} 
\end{definition}
The polynomial $g(z)$ is called the Goppa polynomial. Since $g(z)$ is irreducible and of degree $r$ over $\F_{q^n}$, $g(z)$ does not have any root in $L$ and the code is called an irreducible Goppa code of degree $r$. In this paper $g(z)$ is
always irreducible of degree $r$ over $\F_{q^n}$.\\

It can be shown, see \cite{Chen}, that if $\alpha$ is any root of the Goppa polynomial $g(z)$ then  $\Gamma(L, g)$ is completely described $\alpha$ and a parity check matrix $\bf{H}(\alpha)$ is given by
\begin{equation}
{\bf H}(\alpha)=\left(\frac{1}{\alpha-\zeta_{0}}\frac{1}{\alpha-\zeta_{1}}\cdots \frac{1}{\alpha-\zeta_{q^n-1}}\right),
\end{equation}
where $L=\F_{q^n}=\{\zeta_i: 0\leq i \leq q^n-1 \}$.  

Next we give the definition of an irreducible Goppa code extended with an overall parity check.
\begin{definition}
Let $\Gamma(L,g)$ be an irreducible Goppa code of length $q^n$. Then the extended code $\overline{\Gamma(L,g)}$ is defined by $\overline{\Gamma(L,g)}=\lbrace(c_{0},c_{1},...,c_{q^{n}}): (c_{0},c_{1},...,c_{q^{n}-1})\in \Gamma(L,g)\hspace{0.2cm} \mbox{and}~ \sum_{i=0}^{q^n}c_{i}=0\rbrace $.
\end{definition}
In this paper we take $n=1$. That is, we consider irreducible and extended irreducible Goppa codes of length $q$ and $q+1$, respectively. 
\subsection{Matrices of a given order in $GL(2,q)$}\label{Matriculation}
Let $A\in GL(2,q)$ be of order $D$. We obtain a characterization of the elements of $GL(2,q)$ based on minimal polynomials, conjugacy classes and the order of each matrix. We focus our attention on elements of $GL(2,q)$ which fit our purpose. Elements of $GL(2,q)$ of a given order turn out to be useful in the enumeration of irreducible and extended irreducible Goppa codes.

It is well known that if the order of a matrix $A\in GL(2,q)$ is $D$ then $D\mid p(q-1)$ or $D\mid (q^{2}-1)$ and that the minimal polynomial $m_A(x)$ of $A$ divides $x^D-1$. Combining these facts, Proposition 4.2.2 in \cite{Basheer} and Lemma 2.1 in \cite{Khowa2} we obtain the following theorem.
\begin{theorem}\label{thm_matrix_types}
Let $A=\left(
\begin{array}{cc}
a&b\\
c&d
\end{array}\right)
\in GL(2,q)$ be of order $D$ and $\F_{q}=\F_{p^{t}}$. Denote the minimal polynomial of $A$ by $m_A(x)$. Then  
\begin{enumerate}
\item If $D=1$, then $m_A(x)=x-1$. Thus $A=\left(
\begin{array}{cc}
1&0\\
0&1
\end{array}\right)$.
\item If $p\mid D$ then $m_A(x)=(x-1)^2$ and $A$ is conjugate with a matrix of the form $\left(
\begin{array}{cc}
1&b\\
0&1
\end{array}\right)$ where $b \in \F^\ast_{q}$.
\item If $(p,D)=1$, $D\mid(q-1)$ and $m_A(x)=(x-1)(x-a)$ for some $a \in \F_{q}-\{0,1\}$ where $A$ is not a multiple of the identity matrix, then $A$ is conjugate with a matrix of the form $\left(
\begin{array}{cc}
a&0\\
0&1
\end{array}\right)$ where $a\in \F_{q}-\{0,1\}$.
\item If $(p,D)=1$ and $D\mid(q+1)$ where $D>2$ and $m_A(x)=x^2-\xi x-\zeta \in \F_{q}$ where $m_A(x)$ is irreducible over $\F_{q}$, then $A$ is conjugate with a matrix of the form $\left(
\begin{array}{cc}
0&1\\
\zeta&\xi
\end{array}\right)$. 
\end{enumerate}
\end{theorem}
\subsection{Equivalence Classes}
\subsubsection{The set $\S$}
An irreducible Goppa code can be defined by any root of its Goppa polynomial. As such the set of all roots of such polynomials is important and we make the following definition.

\begin{definition}
The set $\S$ is the set of all elements in $\F_{q^{r}}$ of degree $r$ over $\F_{q} $.  
\end{definition}

\subsubsection{Maps on $\S$}\label{maps_on_S}
 We define the following maps on $\S$.
\begin{definition}\label{Maps_on_s}
 Let $\alpha \in \S$. Mappings of $\alpha$ of Types 1, 2 and 3 are defined as follows:
\begin{itemize}
\item[Type 1] $\sigma^i:\alpha \mapsto \alpha^{q^i}$ where $\sigma$ denotes the Frobenius
automorphism of $\F_{q^{r}}$ leaving $\F_q$ fixed and $0 \leq i\leq r$;\label{uyu1}
\item[Type 2] $\pi_A:\alpha\mapsto a\alpha+b$ where $A=\left(
\begin{array}{cc}
a&b\\
0&1
\end{array}\right)
\in GL(2,q)$.\label{uyu2}
\item[Type 3] $\pi_{B}:\alpha \mapsto \frac{a\alpha+b}{c\alpha+d}$, where $B=\left(
\begin{array}{cc}
a&b\\
c&d
\end{array}\right)
\in GL(2,q)$.\label{uyu3} 
\end{itemize} 
\end{definition}

It has been shown in \cite{Berger} that the composition of Type 1 and Type 2 sends irreducible Goppa codes into equivalent irreducible Goppa codes and the composition of Type 1 and Type 3 maps sends extended irreducible Goppa codes into equivalent extended irreducible Goppa codes. Note that the ``action" of Type 1 and Type 2 on $\S$ was used in \cite{Chen}, Theorem 1, to find bounds on the number of equivalence classes of irreducible Goppa codes.  
\subsubsection{Groups arising from Type 1, Type 2 and Type 3 maps}
In this section we define groups which arise from Type 1, Type 2 and Type 3 maps. The action of these groups on $\S$ will help in counting irreducible and extended irreducible Goppa codes. 
\begin{definition}
Let $G$ denote the set of all maps $\{\sigma^{i}:1\leq i\leq r\}$. $G$ forms a group under the composition of mappings. It is the group of Frobenius automorphisms. It is shown in \cite{Ryan3} that $G$ acts on $\S$.
\end{definition}

\begin{definition}
Let $F$ denote the set of all maps $\left\lbrace\pi_{A}:A=\left(
\begin{array}{cc}
a&b\\
0&1
\end{array}\right)
\in GL(2,q) \right\rbrace$. $F$ forms a group under the composition of mappings and is isomorphic to the group of affine linear transformations.
\end{definition}

Observe that there is an action of the projective linear group $PGL(2,q)$ on $\S$ via the map $\pi_B(\alpha)=[B](\alpha)=\frac{a\alpha+b}{c\alpha+d}$ where $\alpha \in \S$ and $[B]\in PGL(2,q)$, see \cite{Khowa2}.

\subsubsection{Actions of $F$, $PGL(2,q)$ and $G$}\label{Actions}
We first consider the action of the affine group $F$ on $\S$. For each $\alpha \in \S$, the action of $F$ on $\S$ induces orbits denoted $A(\alpha)$ where $A(\alpha)=\{a\alpha+b:a \neq 0, b \in \F_{q}\}$, and called the affine set containing $\alpha$. We denote the set of all affine sets, $ \{A(\alpha):\alpha \in \S \}$, by $\A$. Since $|A(\alpha)|=q(q-1)$ then $|\A|=|\S|/q(q-1)$. It can be shown that $G$ acts on the set $\A$, see \cite{Ryan2}. We will then consider the action of $G$ on $\A$ to obtain orbits in $\S$ of $FG$. The number of orbits in $\S$ under $FG$ will give us an upper bound on the number of irreducible Goppa codes.

Next we consider the action of $PGL(2,q)$ on $\S$. The action of $PGL(2,q)$ on $\S$ induces orbits denoted by $O(\alpha)$ where $O(\alpha)=\{\frac{a\alpha +b}{c\alpha +d}: a,b,c,d\in \mathbb{F}_{q}, ad-bc \neq 0\}$. We will refer to $O(\alpha)$ as a projective linear set. By Theorem 2.3 in \cite{Khowa2}, $|O(\alpha)|=q^{3}-q$.

We denote the set of all projective linear sets in $\mathbb{S}$ under the action of $PGL(2,q)$ by {\rm $\mathbb{O}$. That is, $\mathbb{O}=\{O(\alpha): \alpha \in \mathbb{S}\}$. Observe that $\mathbb{O}$ partitions the set $\mathbb{S}$} and that $G$ acts on the set $\mathbb{O}$ \cite{Ryan3}.

It is shown in \cite{Ryan3} that each projective linear set $O(\alpha)$ in $\mathbb{O}$ can be partitioned into $q+1$ affine sets. See the theorem below.

\begin{theorem}\label{O}
For $\alpha \in \mathbb{S}, O(\alpha)=A(\alpha)\cup A(\frac{1}{\alpha})\cup A(\frac{1}{\alpha + 1})\cup A(\frac{1}{\alpha + \xi _{1}}) \cup A(\frac{1}{\alpha + \xi _{2}})\cup \dots \cup A(\frac{1}{\alpha + \xi _{q-2}})$ where $\mathbb{F}_{q}=\{0,1,\xi_{1}, \xi_{2},\ldots, \xi_{q-2}\}$.
\end{theorem}

Observe that the sets $\mathbb{O}$ and $\mathbb{A}$ are different. $\O$ and $\mathbb{A}$ are both partitions of $\mathbb{S}$ but $|\mathbb{A}|=(q+1) \times |\mathbb{O}|$. 

We will use the actions of $PGL(2,q)$ and $G$ on $\S$ to find an upper bound on the number of extended irreducible Goppa codes. Firstly, we will apply the action of $PGL(2,q)$ on $\S$ to obtain projective linear sets $O(\alpha)$. Then we will consider the action of $G$ on $\O$. The number of orbits in $\O$ under the action of $G$ will give an upper bound on the  number of extended irreducible Goppa codes. To find the number of orbits we will use the Cauchy-Frobenius counting theorem, see \cite{Isaacs}.

The group $G=\langle \sigma\rangle$ is cyclic of order $r$. In analysing the action of $G$ we will make use of the fact that subgroups of $G$ are of the form $\langle \sigma^{r_1}\rangle$ where $r_1\mid r$. Clearly, $|\langle \sigma^{r_1}\rangle|=\frac{r}{r_1}=\bar{r}_1$.

Next we note that if $r=3$ then the number of projective linear sets in $\O$ is $|\O|=\frac{|\S|}{q^{3}-q}=\frac{q^{3}-q}{q^{3}-q}=1$. That is, there is just one projective linear set containing all $q^{3}-q$ elements of $\S$. We put the result in a theorem.
\begin{theorem}
When $r=3$ the set $\S$ consists of just one projective linear set, that is, $\S = O(\alpha)$ for any $\alpha \in \S$.
\end{theorem}

\begin{corollary}\label{Coro_r=3}
All extended irreducible Goppa codes of length $q+1$ with $r=3$ are equivalent.
\end{corollary}
The result in Corollary \ref{Coro_r=3} is well known, for example see \cite{Chen} and \cite{Moreno}. 
 
Now, suppose that $[B](\alpha)=\frac{a\alpha+b}{c\alpha+d}=\alpha^{q^{r_1}}$ where $r_1\mid r$. We see that if $D$ is the smallest positive integer such that $\alpha^{q^{Dr_1}}=\alpha$, then $D\mid r$ since $\alpha \in \F_{q^{r}}$. Observe that $\alpha^{q^{Dr_1}}=\alpha$ implies that $\alpha^{q^{Dr_1}}=[B^D](\alpha)=[I_2](\alpha)=\alpha$. Thus the order $D$ of $B\in GL(2,q)$ is $D=\frac{r}{r_1}=\bar{r}_1$.

\section{Enumeration of irreducible Goppa codes}
We count irreducible Goppa codes by using the tools developed in \cite{Ryan2} where an upper bound on the number of irreducible Goppa codes of degree $r$ and length $q^n$ is given. The upper bound is found by counting the number of affine sets in $\mathbb{A}$ fixed under the action of subgroups $G$ and then applying the Cauchy-Frobenius Theorem. As opposed to \cite{Ryan2}, where the upper bound is given in the form of an algorithm, we obtain analytic formulas for the upper bound. 

Suppose $A(\alpha)$ is fixed by $\langle\sigma^{r_1}\rangle$. Then $\sigma ^{r_1}(A(\alpha))=A(\alpha)$. So we have $\sigma^{r_1}(\alpha)=\alpha ^{q^{r_1}}=\zeta \alpha + \xi$ for some $\zeta\neq 0,\xi \in \mathbb{F}_{q}$. Thus $\alpha^{q^{r_1}}=[A](\alpha)$ where $A=\left(
\begin{array}{cc}
\zeta&\xi\\
0&1
\end{array}\right)\in GL(2,q)$. Now, $\alpha^{q^{{r_1}\bar{r}_1}}=[A^{\bar{r}_1}](\alpha)=[I_2](\alpha)=\alpha$ so $A$ is a matrix of order $\bar{r}_1$. We divide our analysis according to whether $\bar{r}_1=1$, $\mbox{gcd}(p,\bar{r}_1)=1$ and $\mbox{gcd}(p,\bar{r}_1)=p$. 

\subsection{Affine sets fixed under $\langle\sigma^{r_1}\rangle$ when $\bar{r}_1=1$}\label{Aff_fixed_D=1}
If $\bar{r}_1=1$ then $r_1=r$ and  $\alpha ^{q^{r_1}}=[A](\alpha)=[I_2](\alpha)=\alpha$ and $A=\left(
\begin{array}{cc}
1&0\\
0&1
\end{array}\right)$, see Theorem \ref{thm_matrix_types}. Now $\alpha ^{q^{r_1}}=\alpha$ if and only if $\F_{q^{r_1}}$ contains elements of $\S$. We know that this is true if and only if $r_1=r$. It is easy to see that every affine set is fixed under $\langle\sigma^{r}\rangle$. By Corollary 3.5 in \cite{Ryan_Disc} the number of affine sets fixed under $\langle\sigma^{r}\rangle$ is 
$$\frac{|\S|}{q^{2}-q}.$$ 
\begin{example}
Let $q=2^6$ and $r=3$. There are 17 affine sets in $\A$ and all of them are fixed under $\langle\sigma^{3}\rangle$. 
\end{example}

\subsection{Affine sets fixed under $\langle\sigma^{r_1}\rangle$ when $\mbox{gcd}(p,\bar{r}_1)=1$}
Now suppose that $A(\alpha)\in \A$ is fixed under $\langle\sigma^{r_1}\rangle$ where $\mbox{gcd}(p,\bar{r}_1)=1$. Then we have that  $\alpha ^{q^{r_1}}=[A](\alpha)$ and by Theorem \ref{thm_matrix_types}, we may take $A=\left(
\begin{array}{cc}
\zeta&0\\
0&1
\end{array}\right)\in GL(2,q)$ of order $\bar{r}_1$. Thus $\alpha^{q^{r_1}}=\zeta\alpha$ and as such $\alpha$ satisfies an equation of the form 
\begin{equation}\label{Fac_Eqn_coprime}
x^{q^{r_1}}-\zeta x=0,
\end{equation}
where $\zeta$ is of order $\bar{r}_1$.  

Next we note that the factorization of $F_{r_1}(x)=x^{q^{r_1}}-\zeta x$ was considered in \cite{Gare}. Using our notation and Theorem 4 in \cite{Gare} we obtain the following result.

\begin{theorem}\label{T}
Let $\bar{r}_1>1$ and suppose $(p,\bar{r}_1)=1$ and that $\bar{r}_1 \mid (q-1)$. Let $r=\bar{r}_1 u, u \in \N$ and $T(r)$ be the set of all roots of irreducible factors of degree $r$ in the factorization of $F_{r_1}(x)=x^{q^{s}-1}-\zeta \in \F_{q}[x]$, where $\zeta\neq 1 \in \F^{\ast}_{q}$ is of order $\bar{r}_1$. If $r \not\equiv 0 \pmod{\bar{r}_1}$ then $|T(r)|=0$. Otherwise  
 $$|T(r)|=\sum_{\substack{d|u\\(d,\bar{r}_1)=1}}\mu(d)\left(q^{\frac{u}{d}}-1\right).$$ 
\end{theorem}

Observe that since $A=\left(
\begin{array}{cc}
\zeta&0\\
0&1
\end{array}\right)\in GL(2,q)$ is of order $\bar{r}_1$ there are $\phi(\bar{r}_1)$ conjugacy classes in $GL(2,q)$ in this case and a polynomial arising from a representative of each conjugacy class contributes $|T(r)|$ roots to $\S$. Thus there are $\phi(\bar{r}_1)|T(r)|$ roots which lie in $\S$. Note that this closed formula is a partial answer to Remark 4.5 in \cite{Ryan_Disc}. 
\begin{example}
Let $q=2^6$ and $r=6$. There are 672 irreducible factors of degree 6 in the factorization of $F_2(x)=x^{2^{12}}-\varepsilon^{42}x$ where $\varepsilon$ is a primitive element of $\F_{2^6}$. Hence there are 8,064 roots of polynomials of the form $F_2(x)$ which lie in $\S$ where $A=\left(
\begin{array}{cc}
\varepsilon^{42}&0\\
0&1
\end{array}\right)\in GL(2,2^6)$ is of order $3$.
\end{example}

Next we find the number of affine sets fixed under $\langle\sigma^{r_1}\rangle$. We know that if $\bar{r}_1>1$ where $p\nmid \bar{r}_1$ and $\bar{r}_1\mid (q-1)$ then there are $\phi(\bar{r}_1)|T(r)|$ roots of the equations of the form $x^{q^{r_1}}-\zeta x=0$ which lie in $\S$. By \cite[Theorem 4.4]{Ryan_Disc}, each polynomial $F_{r_1}(x)$ has $q-1$ roots in exactly one $A(\alpha)$. Thus we have proved the following theorem. 

\begin{theorem}\label{Aff_sets_fixed_D_divides}
Suppose $\bar{r}_1>1$ where $p\nmid \bar{r}_1$ and $\bar{r}_1\mid (q-1)$. Then there are $\tau=\frac{\phi(\bar{r}_1)|T(r)|}{q-1}$ affine sets fixed by $\langle\sigma^{r_1}\rangle$.
\end{theorem}

\subsection{Affine sets fixed under $\langle\sigma^{r_1}\rangle$ when $\mbox{gcd}(p,\bar{r}_1)=p$}
Suppose that $A(\alpha)\in \A$ is fixed under $\langle\sigma^{r_1}\rangle$ where $\mbox{gcd}(p,\bar{r}_1)=p$. Then we have that  $\alpha ^{q^{r_1}}=[A](\alpha)$ and by Theorem \ref{thm_matrix_types}, we may take $A=\left(
\begin{array}{cc}
1&\beta\\
0&1
\end{array}\right)\in GL(2,q)$ where $\beta\in \F^{\ast}_q$. Thus $\alpha^{q^{r_1}}=\alpha+\beta$ and as such $\alpha$ satisfies an equation of the form 
\begin{equation}\label{Fac_Eqn_p}
x^{q^{r_1}}-x-\beta =0.
\end{equation}

Observe that if $\bar{r}_1=p$ hence $r_1=\frac{r}{p}$ then we will take $A=\left(
\begin{array}{cc}
1&1\\
0&1
\end{array}\right)
$ since we can show, by direct computation, that the order of $A$ is $p$. As such, we have $\alpha ^{q^{\frac{r}{p}}}=\alpha + 1$ and we may assume that $\alpha$ satisfies an equation of type 

\begin{equation}\label{Fac_Eqn_5}
x^{q^{\frac{r}{p}}}-x-1=0.
\end{equation}

Next we note that the factorization of polynomials of the form $F_{r_1}(x)=x^{q^{r_1}}-x-1$ was considered in \cite{Gare}. Using our notation and \cite[Theorem 2]{Gare} we obtain the following result.

\begin{theorem}\label{V}
Suppose $F_{r_1}(x)=x^{q^{r_1}}-x-1$ where $r_1=\frac{r}{p}$ and $\bar{r}_1=p$. Let $r=pu, u \in \N$, and $U(r)$ be the set of roots of irreducible factors of degree $r$ in the factorization of $F_{r_1}(x)$. If $r \not\equiv 0 \pmod{p}$, then $|U(r)|=0$. Otherwise
$$|U(r)|=\sum_{\substack{d|u\\d \not\equiv 0\pmod{p}}}\mu(d)q^{u/d}.$$ 
\end{theorem}
It is easy to see that if $\alpha$ satisfies Equation \ref{Fac_Eqn_5} then all the $q$ elements of the set $\{\alpha +\xi: \xi \in \F_q\}$ also satisfy \ref{Fac_Eqn_5} while the remaining elements in $A(\alpha)$ do not. Hence if $\alpha$ satisfies equation \ref{Fac_Eqn_5} then $A(\alpha)$ contains precisely $q$ roots of Equation \ref{Fac_Eqn_5}. We have proved the following theorem.
\begin{theorem}\label{V1}
If $\bar{r}_1=p$ then there are $\frac{|U(r)|}{q}$ affine sets fixed by $\langle\sigma^{r_1}\rangle$ where $r_1=\frac{r}{p}$.  
\end{theorem}

\begin{corollary}
If $r=p$ and $\bar{r}_1=p$ then $r_1=1$ and there is one affine set fixed by $\langle\sigma\rangle$.
\end{corollary}

\begin{example}
Let $q=2^6$ and $r=6$. Then $|U(r)|=262,080$ and there are $\frac{262,080}{64}=4,095$ affine sets fixed by $\langle\sigma^{3}\rangle$.
\end{example}

Now suppose that $\bar{r}_1=p^2$. We consider $r_1=\frac{r}{p^2}$. Thus, we have $\alpha^{q^{r_1}}=[A](\alpha)$. So $\alpha ^{q^{r}}=[A^{p^2}](\alpha)=[I_2](\alpha)=\alpha$ and so the order of $A$ divides $p^2$. We know that matrices of order $p^2$ do not exist. So we consider matrices of order $p$. We obtain an equation of the form $x^{q^{\frac{r}{p^2}}}-x-1=0$ and all roots of this equation lie in $\F_{q^{\frac{r}{p}}}-\F_{{q^{\frac{r}{p^2}}}}$ and not in $\S$. We have the following theorem.

\begin{theorem}
If $\bar{r}_1=p^2$ then $p^2 \mid r$ and there is no polynomial of degree $r$ in the factorization of $F_{r_1}(x)=x^{q^{r_1}}-x-1$ where $r_1=\frac{r}{p^2}$.
\end{theorem}

Next suppose that $\bar{r}_1=pp_1$ where $p_1$ is some other divisor of $r$. Then $\alpha^{q^{r_1}}=[A](\alpha)$ and $\alpha^{q^{r_1pp_1}}=[A^{pp_1}](\alpha)=[I_2](\alpha)=\alpha$. Then $A^{pp_1}=I_2$. If we take $B=A^p$ of order $p_1$ then we have $\alpha^{q^{sp_1}}=[B^{p_1}](\alpha)$ then $\alpha$ satisfies Equation \ref{Fac_Eqn_coprime}. Also if we take $B=A^{p_1}$ of order $p$ then we have $\alpha^{q^{r_1p}}=[B^{p}](\alpha)$ then $\alpha$ satisfies Equation \ref{Fac_Eqn_5} and this is not possible. So there is no irreducible polynomial of degree $r$ in the factorization of $F_{r_1}(x)$ in this case. We have the following result.
\begin{theorem}
There is no affine set fixed under $\langle\sigma^{r_1}\rangle$ if $\bar{r}_1=p^2$ and $p^2 \mid r$ or $\bar{r}_1=pp_1$ where $p_1$ is some other divisor of $r$.
\end{theorem}
\begin{example}
Let $q=2^4$ and $r=12$. The subgroups $\langle\sigma^{2}\rangle$ and $\langle\sigma^{3}\rangle$ where $\bar{r}_1=6$ and $\bar{r}_1=4$ respectively do not fix any affine set. 
\end{example}
Putting the results together, we have proved the following:
\begin{theorem}\label{main_theo_Gamba_Aff}
With the notation we have established: 
\begin{enumerate}
\item There are $\frac{|\S|}{q^{2}-q}$ affine sets fixed by $\langle\sigma^{r}\rangle$.
\item There are $\frac{\phi(\bar{r}_1)|T(r)|}{q-1}$ affine sets fixed by $\langle\sigma^{r_1}\rangle$  if $\bar{r}_1\mid(q-1)$. 
\item If $(\bar{r}_1,r)=p$ then 
\begin{enumerate}
\item there are $\frac{|U(r)|}{q}$ affine sets fixed by $\langle\sigma^{r_1}\rangle$ if $\bar{r}_1=p$.
\item there is $1$ affine set fixed by $\langle\sigma^{r_1}\rangle$ if $\bar{r}_1=r$.
\end{enumerate}
\end{enumerate}
\end{theorem}
\begin{remark}
This result agrees with \cite[Theorem 4.13]{Ryan_Disc} for $n=1$. Our main contribution here is that we have found closed formulas for $|T(r)|$ and $|U(r)|$.
\end{remark}
\section{Counting extended irreducible Goppa codes}
\subsection{Strategy for counting extended irreducible Goppa codes}
We will use the actions of $PGL(2,q)$ and $G$ on $\S$ to find the maximum number of extended irreducible Goppa codes. Firstly, we will apply the action of the group $PGL(2,q)$ on $\S$ to obtain projective linear sets $O(\alpha)$. Then we will consider the action of $G$ on $\O$. The number of orbits in $\O$ under the action of $G$ will give us an upper bound on the number of extended Goppa codes. 

Recall that a projective linear set can be decomposed as $O(\alpha)=A(\alpha)\cup A(\frac{1}{\alpha})\cup A(\frac{1}{\alpha + 1})\cup A(\frac{1}{\alpha + \xi _{1}}) \cup A(\frac{1}{\alpha + \xi _{2}})\cup \dots \cup A(\frac{1}{\alpha + \xi _{q-2}})$ where $\mathbb{F}_{q^{n}}=\{0,1,\xi_{1}, \xi_{2},\ldots, \xi_{q-2}\}$. Observe that if a projective linear set $O(\alpha)\in \O$ is fixed under $\langle \sigma^{r_1}\rangle$ then $\langle \sigma^{r_1}\rangle$ acts on $O(\alpha)=A(\alpha)\cup A(\frac{1}{\alpha})\cup A(\frac{1}{\alpha + 1})\cup A(\frac{1}{\alpha + \xi _{1}}) \cup A(\frac{1}{\alpha + \xi _{2}})\cup \dots \cup A(\frac{1}{\alpha + \xi _{q-2}})$ and partitions this set of $q+1$ affine sets. We see that some projective linear sets fixed under $\langle \sigma^{r_1}\rangle$ contain fixed affine sets and there is also a possibility of having a fixed projective linear set that does not contain fixed affine sets. We will consider the following possibilities: $\bar{r}_1=1$; $\mbox{gcd}(p,\bar{r}_1)=1$ and $\bar{r}_1\mid (q-1)$; $\mbox{gcd}(p,\bar{r}_1)=p$; and $\mbox{gcd}(p,\bar{r}_1)=1$ where $\bar{r}_1\mid (q+1)$. We will discuss each of the four cases in separate sections.
\subsection{Projective linear sets fixed when $D=1$}\label{trivial}
If $\bar{r}_1=1$ then $\alpha ^{q^{r_1}}=[B](\alpha)=[I_2](\alpha)=\alpha$ and $B=\left(
\begin{array}{cc}
1&0\\
0&1
\end{array}\right)$, see Theorem \ref{thm_matrix_types}. Now $\alpha ^{q^{r_1}}=\alpha$ if and only if $\F_{q^{r_1}}$ contains elements of $\S$. That is, $r_1=r$. By Section \ref{Aff_fixed_D=1} we know that the number of affine sets fixed under $\langle\sigma^{r}\rangle$ where $\bar{r}_1=1$ is 
$\frac{|\S|}{q^{2}-q}$. 

By an argument similar to the one in \cite[Section 4.3.2]{Khowa2}, we find that the number of projective linear sets $O(\alpha)$ fixed under $\langle\sigma^{r}\rangle$ is $\frac{|{\mathbb S}|}{q^3-q}$. 

\subsection{Projective linear sets fixed when $\mbox{gcd}(p,\bar{r}_1)=1$ and $\bar{r}_1\mid (q-1)$ }\label{mid}
Suppose $\mbox{gcd}(p,\bar{r}_1)=1$ and that $\bar{r}_1\mid (q-1)$. Theorem \ref{Aff_sets_fixed_D_divides} gives the number of affine sets fixed by $\langle\sigma^{r_1}\rangle$ in this case.

Next we find the number of projective linear sets $O(\alpha)$ fixed under $\langle\sigma^{r_1}\rangle$. We will do this by finding how many affine sets fixed under $\langle\sigma^{r_1}\rangle$ lie in each fixed projective linear set. 

Suppose $O(\alpha)\in {\mathbb O}$ is fixed under $\langle\sigma^{r_1}\rangle$. Then $\langle\sigma^{r_1}\rangle$ acts on $O(\alpha)=A(\alpha)\cup A(\frac{1}{\alpha})\cup A(\frac{1}{\alpha+1})\cup A(\frac{1}{\alpha+\xi_{1}}) \cup A(\frac{1}{\alpha+\xi_{2}}) \cup A(\frac{1}{\alpha+\xi_{3}})\cup\dots\cup A(\frac{1}{\alpha+\xi_{q-2}})$, a set of $q+1$ affine sets. $\langle\sigma^{r_1}\rangle$ partitions this set of $q+1$ affine sets. The possible lengths of an orbit are $1$ and factors of $\bar{r}_1$. Now, since $\bar{r}_1\mid q-1$ then $q+1=q-1+2\equiv 2 \pmod{\bar{r}_1}$. We claim that each $O(\alpha)$ fixed under $\langle\sigma^{r_1}\rangle$ contains 2 affine sets which are fixed under $\langle\sigma^{r_1}\rangle$. Observe that if $\bar{r}_1$ is prime then we are done. So we will suppose that $\bar{r}_1$ is composite. That is, we can find non-negative integers $e_1,e_2,\dots,e_t$ such that $d_1e_1+d_2e_2+\cdots+d_te_t=q+1$ where $d_i\mid \bar{r}_1$, $1\leq i\leq t$, $d_1=1$ and $d_t=\bar{r}_1$. Note that we can always choose a factor $d_i$, $1< i<t$ such that $\langle\sigma^{d_ir_1}\rangle$ is of prime order $\frac{\bar{r}_1}{d_i}$ and fixes $O(\alpha)$. Now, an $O(\alpha)$ fixed under $\langle\sigma^{d_ir_1}\rangle$ contains two fixed affine sets so it follows that an $O(\alpha)$ fixed under $\langle\sigma^{r_1}\rangle$ also contains two fixed affine sets. Thus if $\bar{r}_1\mid (q-1)$ then each $O(\alpha)$ fixed under $\langle\sigma^{r_1}\rangle$ contains 2 fixed affine sets.

We have the following theorem.
\begin{theorem}\label{thm_tau_ext}
Let $\mbox{gcd}(p,\bar{r}_1)=1$ and $\bar{r}_1\mid (q-1)$. The number of projective linear sets fixed under $\langle\sigma^{r_1}\rangle$ is $\frac{\tau}{2}$ where $\tau$ is defined in Theorem \ref{Aff_sets_fixed_D_divides}. 
\end{theorem}

\subsection{Projective linear sets fixed when $\mbox{gcd}(p,\bar{r}_1)=p$}\label{nmid}
In this section we obtain the number of projective linear sets fixed when $\mbox{gcd}(p,\bar{r}_1)=p$. By Theorem \ref{V1}, there are $\frac{|U(r)|}{q}$ affine sets fixed by $\langle\sigma^{r_1}\rangle$ when $\mbox{gcd}(p,\bar{r}_1)=p$. We will do this by finding how many affine sets fixed under $\langle\sigma^{r_1}\rangle$ lie in a projective linear set fixed under $\langle\sigma^{r_1}\rangle$. 

We claim that each of the $O(\alpha)$ in $\O$ fixed under $\langle\sigma^{r_1}\rangle$ contains precisely one affine set which is fixed under $\langle\sigma^{r_1}\rangle$. It suffices to show that $O(\alpha)$ cannot contain two affine sets which are fixed under $\langle\sigma^{r_1}\rangle$. Without loss of generality, suppose $A(\alpha)$ is fixed under $\langle\sigma^{r_1}\rangle$. Recall that $O(\alpha)=A(\alpha)\cup A(\frac{1}{\alpha})\cup A(\frac{1}{\alpha+1})\cup A(\frac{1}{\alpha+\xi_{1}}) \cup A(\frac{1}{\alpha+\xi_{2}}) \cup A(\frac{1}{\alpha+\xi_{3}})\cup\dots\cup A(\frac{1}{\alpha+\xi_{q-2}})$. We show that none of the affine sets after $A(\alpha)$ in the above decomposition of $O(\alpha)$ is fixed under $\langle\sigma^{r_1}\rangle$. This is done by showing that no element in any of these affine sets satisfies Equation \ref{Fac_Eqn_5}. It is sufficient to show that no element in $A(\frac{1}{\alpha})$ satisfies $x^{q^{\frac{r}{p}}}-x-1=0=0$. A typical element of $A(\frac{1}{\alpha})$ has the form $\frac{\zeta}{\alpha}+ \xi$ and substituting this into $x^{q^{\frac{r}{p}}}-x-1$ we get $(\frac{\zeta}{\alpha}+ \xi)^{q^{\frac{r}{p}}}-(\frac{\zeta}{\alpha}+ \xi)-1=\frac{-\alpha^{2}-\alpha- \zeta}{\alpha^{2}+\alpha}\neq 0$, since $\alpha$ is an element of degree $r>3$ over $\mathbb{F}_{q}$. We conclude that $A(\frac{1}{\alpha})$ is not fixed under $\langle\sigma^{r_1}\rangle$ and in fact $A(\alpha)$ is the only affine set in $O(\alpha)$ fixed under $\langle\sigma^{r_1}\rangle$. It follows that the number of projective linear sets $O(\alpha)$ in $\O$ which are fixed under $\langle\sigma^{r_1}\rangle$ where $(p,\bar{r}_1)=p$ is $\frac{|U(r)|}{q}$. Thus we have proved the following.
\begin{theorem}
If $\mbox{gcd}(p,\bar{r}_1)=p$ and $\bar{r}_1=p$, then the number of projective linear sets fixed under $\langle\sigma^{r_1}\rangle$ is the same as the number of affine sets fixed under $\langle\sigma^{r_1}\rangle$. 
\end{theorem}
 
\subsection{Projective linear sets fixed by $\langle\sigma^{r_1}\rangle$ where $\mbox{gcd}(p,\bar{r}_1)=1$ and $\bar{r}_1\mid (q+1)$}\label{no_aff_fixed}
Suppose that $O(\alpha)\in \O$ is fixed by $\langle\sigma^{r_1}\rangle$ where $\mbox{gcd}(p,\bar{r}_1)=1$, $\bar{r}_1\mid (q+1)$ and $\bar{r}_1>2$. Then we have that $\alpha ^{q^{r_1}}=[A](\alpha)$, where $\mbox{gcd}(p,\bar{r}_1)=1$ and $\bar{r}_1\mid (q+1)$. Then, by Theorem \ref{thm_matrix_types}, $A$ is conjugate with a matrix of the form $B=\left(
\begin{array}{cc}
0&1\\
\zeta&\xi
\end{array}\right)\in GL(2,q)
$ where the minimal polynomial of $B$, $m_B(x)$, is an irreducible quadratic polynomial over $\F_{q}$. Without loss of generality, we will take $A=\left(
\begin{array}{cc}
a&b\\
c&d
\end{array}\right)\in GL(2,q)$ where, as above, $m_A(x)$ is an irreducible quadratic polynomial over $\F_{q}$.

Now, $\alpha^{q^{r_1}}=[A](\alpha)=\frac{a\alpha+b}{c\alpha+d}$ implies that $\alpha$ satisfies an equation of the form 
\begin{equation}\label{eqnyowina}
F_{r_1}(x)=cx^{q^{r_1}+1}+dx^{q^{r_1}}-ax-b=0.
\end{equation} It is clear from the foregoing discussion that in order to find the number of projective linear sets fixed under $\langle \sigma^{r_1}\rangle$ we need to find roots of Equation \ref{eqnyowina} which lie in $\S$. Observe that in this case there is no affine set fixed in the decomposition of $O(\alpha)$.

Note that there are $\frac{\phi(\bar{r}_1)}{2}$ polynomials of the form $F_{r_1}(x)=cx^{q^{r_1}+1}+dx^{q^{r_1}}-ax-b\in \F_{q}[x]$ each of which corresponds to a representative of a conjugacy class of matrices of order $\bar{r}_1$ where the minimal polynomials of such matrices are irreducible quadratic polynomials over $\F_{q}$, see Theorem 2.2 in \cite{Khowa2}.

We now consider the factorization of $F_{r_1}(x)=cx^{q^{r_1}+1}+dx^{q^{r_1}}-ax-b\in \F_{q}[x]$. We begin by considering the factorization of $F_{r_1}(x)$ where $\mbox{gcd}(p,\bar{r}_1)=1$, $\bar{r}_1\mid (q+1)$ and $\bar{r}_1>2$ is even. Note that the assumption that $\bar{r}_1\mid (q+1)$ where $\bar{r}_1$ is even implies that the characteristic of $\F_q$ is odd.

Suppose that $\alpha$ is a root of $F_{r_1}(x)$ where $\bar{r}_1>2$ is even. Then $\bar{r}_1=2d$ where $d>1$ is an integer. Thus we have $\alpha ^{q^{r_1\bar{r}_1}}=[A^{\bar{r}_1}](\alpha)=[I_2](\alpha)=\alpha$. So $A^{2d}=(A^d)^2=B^2=I_2$, where $B=A^d$. Now, since $B^2=I_2$, without loss of generality we can take $B=\left(
\begin{array}{cc}
q-1&0\\
0&q-1
\end{array}\right)
$ since the only elements $\zeta \in \F_{q}$ such that $\zeta^2=1$ are $\zeta=1$ and $\zeta=q-1$. Thus $\alpha$ satisfies an equation of the form $F_{r_1}(x)=(q-1)(x^{q^{r_1}}-x)=0$. By the argument in Section \ref{trivial} there are no irreducible polynomials of degree $r$ in the factorization of $F_{r_1}(x)=(q-1)(x^{q^{r_1}}-x)$. We have proved the following.
\begin{theorem}
Suppose $\mbox{gcd}(p,\bar{r}_1)=1$, $\bar{r}_1\mid (q+1)$ and $\bar{r}_1=2d$ where $d>1$ is an integer. Then there is no projective linear set fixed under $\langle\sigma^{r_1}\rangle$.    
\end{theorem}

 \begin{example}
Consider $q=3^3$ and $r=4$. There are no polynomials of degree $4$ in the factorization of $F_{1}(x)=2x^{28}+ 2\in \F_{3^3}[x]$ where $A=\left(
\begin{array}{cc}
0&1\\
2&0
\end{array}\right)
\in GL(2,3^3)$ is of order $4$. Hence there is no projective linear set fixed under $\langle\sigma\rangle$. 
\end{example}
Next we consider the factorization of $F_{r_1}(x)$ where $\mbox{gcd}(p,\bar{r}_1)=1$, $\bar{r}_1\mid (q+1)$ and $d\mid \bar{r}_1$ but $d^2\nmid \bar{r}_1$.
Suppose that $\alpha$ is a root of $F_{r_1}(x)$ and that $d\mid \bar{r}_1$ where $1<d<\bar{r}_1$. We have $\alpha ^{q^{r_1}}=[A](\alpha)$ where $A$ is a matrix of order $\bar{r}_1$. Now, if $d\mid \bar{r}_1$ then $\alpha ^{q^{r_1d}}=[A^d](\alpha)=[B](\alpha)$, where $B=A^d$ and the order of $B$ is $\frac{\bar{r}_1}{d}$ since $\mbox{GCD}(\bar{r}_1,d)=d$. We see that for $\alpha$ to satisfy $F_{r_1d}(x)$, $d$ must divide $\frac{\bar{r}_1}{d}$. That is $\alpha$ satisfies $F_{r_1d}(x)$ provided $d^2\mid \bar{r}_1$. We have proved the following theorem.
\begin{theorem}
Suppose $\alpha$ satisfies $F_{r_1}(x)$ and $d\mid \bar{r}_1$ where $1<d<\bar{r}_1$. If $d^2\nmid \bar{r}_1$ then $\alpha$ does not satisfy $F_{r_1d}(x)$.
\end{theorem}
This theorem is significant because it tells us about the existence of an irreducible factor of degree $r$ in the factorization of $F_{r_1}(x)$ by looking at the factorization of $F_{r_1d}(x)$.

Now we consider when an irreducible quadratic polynomial of the form $cx^2+(d-a)x-b\in \F_{q}[x]$ divides $F_{r_1}(x)$. 
\begin{proposition}\label{quad_factor}
 If $P(x)=cx^2+(d-a)x-b\in \F_{q}[x]$ is irreducible, then $P(x)\mid F_{r_1}(x)$ if and only if $r_1=2w$.
\begin{proof}
Suppose $P(x)=cx^2+(d-a)x-b\in \F_{q}[x]$ is irreducible and that $P(\alpha)=0$. Then $\alpha^{q^{2}}=\alpha$. If $r_1$ is even then $\mbox{gcd}(2,r_1)=2$ and this implies that $\alpha^{q^{r_1}}=\alpha$. Thus $F_{r_1}(\alpha)=c\alpha^{q^{r_1}+1}+d\alpha^{q^{r_1}}-a\alpha-b=c\alpha^2+d\alpha-a\alpha-b=0$. So $P(x)\mid F_{r_1}(x)$.
% 
%  Now, if $s$ is odd and $P(x)=cx^2+(d-a)x-b$ divides $F_s(x)$ then it follows that $P(x)$ is reducible since no root of $P(x)$ lies in $\F_{q^s}$ in this case.\\
% 
If $r_1$ is odd and $P(\alpha)=0$ then $\alpha^{q^{r_1}+1}=\alpha^{q+1}$. So $F_{r_1}(\alpha)=c\alpha^{q^{r_1}+1}+d\alpha^{q^{r_1}}-a\alpha-b=c\alpha^{q+1}+d\alpha^{q}-a\alpha-b=(\alpha^{q-1}-1)(c\alpha^2+d\alpha)+c\alpha^2+d\alpha-a\alpha-b=0$. We know that $c\alpha^2+d\alpha-a\alpha-b=0$ so $F_{r_1}(\alpha)=0\Leftrightarrow \alpha^{q-1}-1=0$. This means $\alpha \in \F_{q}$ and $P(x)$ is reducible.    
\end{proof} 
\end{proposition}
 
Consequently, the parity of $r_1$ indicates whether or not there is a quadratic factor in the factorization of $F_{r_1}(x)$. The following theorem, see \cite{Reis} for its proof, gives the factorization of $F_{r_1}(x)=cx^{q^{r_1}+1}+dx^{q^{r_1}}-ax-b\in \F_{q}[x]$.
        
\begin{theorem}\label{main_theo_1}
Let $\mbox{gcd}(p,\bar{r}_1)=1$ and $\bar{r}_1\mid (q+1)$ where $\bar{r}_1\neq 2t$, $t>1$. Let $X(r)$ be the set of roots of polynomials of degree $r=\bar{r}_1  u$ in the factorization of $F_{r_1}(x)=cx^{q^{r_1}+1}+dx^{q^{r_1}}-ax-b$ where the minimal polynomial of $A$ is an irreducible quadratic polynomial over $\F_{q}$. Then
\[|X(r)|=\sum_{\substack{d|u\\d \not\equiv 0\pmod{\bar{r}_1  }}}\mu(d)(q^{\frac{u}{d}}-(-1)^{\frac{u}{d}}).\] 
\end{theorem}
\begin{example}
Consider $q=2^5$ and $r=6$. We want to find the number of irreducible polynomials of degree 6 in the factorization of $F_{2}(x)=x^{2^{10}+1}-x+1\in \F_{2^5}[x]$. Using Theorem \ref{main_theo_1}, with $r=\bar{r}_1 u=3\times 2$ there are $\frac{990}{6}=165$ polynomials of degree 6 in the factorization of $F_{2}(x)=x^{2^{10}+1}+x+1\in \F_{2^5}[x]$.
\end{example}

Theorem \ref{main_theo_1} implies that each polynomial $F_{r_1}(x)=cx^{q^{r_1}+1}+dx^{q^{r_1}}-ax-b$ contributes $|X(r)|$ roots to $\S$. Recall that there are $\frac{\phi(\bar{r_1})}{2}$ conjugacy classes of matrices of order $\bar{r_1}$ whose eigenvalues lie in $\F_{q^{2}}$ hence there are $\frac{\phi(\bar{r_1})}{2}$ polynomials counting representatives only. If we let $\S_F$ be the set of roots of the $\frac{\phi(\bar{r_1})}{2}$ polynomials $F_{r_1}(x)$ which lie in $\S$ then $|\S_F|=\frac{\phi(\bar{r_1})|X(r)|}{2}$.

In the following theorem we count the number of roots of $F_{r_1}(x)$ which lie in $O(\alpha)$.  
\begin{theorem}\label{thm_mbambande}
If $O(\alpha)$ is a projective linear set fixed by $\langle\sigma^{r_1}\rangle$, then $O(\alpha)$ contains $q+1$ roots of $F_{r_1}(x)=cx^{q^{r_1}+1}+dx^{q^{r_1}}-ax-b$.
\begin{proof}
Suppose $O(\alpha)$ is fixed under $\langle\sigma^{r_1}\rangle$, then $\alpha^{q^{r_1}}=[A](\alpha)=\frac{a\alpha+b}{c\alpha+d}$, where $A=\left(
\begin{array}{cc}
a&b\\
c&d
\end{array}\right) \in GL(2,q)
$.\\ Thus we may assume that $\alpha$ satisfies $F_s(x)=cx^{q^{r_1}+1}+dx^{q^{r_1}}-ax-b$.\\ Recall that $O(\alpha)=A(\alpha)\cup A(\frac{1}{\alpha})\cup A(\frac{1}{\alpha + 1})\cup A(\frac{1}{\alpha + \nu _{1}}) \cup A(\frac{1}{\alpha + \nu _{2}})\cup \dots \cup A(\frac{1}{\alpha + \nu _{q-2}})$ where $\mathbb{F}_{q}=\{0,1,\nu_{1}, \nu_{2},\cdots, \nu_{q-2}\}$.\\ It is easy to see that $\alpha$ is the only element in $A(\alpha)$ which satisfies $F_{r_1}(x)$. Now, let us fix $\nu \in \F_{q}$, we would like to count how many elements of the affine set $A(\frac{1}{\alpha + \nu})$ satisfy $F_{r_1}(x)$. A typical element of $A(\frac{1}{\alpha + \nu})$ is $\frac{\zeta}{\alpha + \nu}+\xi$ where $\zeta\neq 0,\xi \in \F_{q}$. If $\frac{\zeta}{\alpha + \nu}+\xi$ is a root of $F_{r_1}(x)=cx^{q^{r_1}+1}+dx^{q^{r_1}}-ax-b$ then
$$c\left( \frac{\zeta}{\alpha + \nu}+\xi\right)^{q^{r_1}+1}+d\left( \frac{\zeta}{\alpha + \nu}+\xi\right)^{q^{r_1}}-a\left( \frac{\zeta}{\alpha + \nu}+\xi\right)-b=0.$$

From this we obtain 
 \begin{eqnarray}\label{array0}
\left( \frac{\zeta}{\alpha + \nu}+\xi\right)^{q^{r_1}}=\frac{a\left( \frac{\zeta}{\alpha + \nu}+\xi\right)+b}{c\left( \frac{\zeta}{\alpha + \nu}+\xi\right)+d}.
\end{eqnarray}
So
\begin{eqnarray}\label{array1}
\frac{\zeta}{\alpha^{q^{r_1}}+\nu}+\xi=\frac{a\zeta+a\xi\alpha+a\xi\nu+b\alpha+b\nu}{c\zeta+c\xi\alpha+c\xi\nu+d\alpha+d\nu}.
\end{eqnarray}
Thus $\nu \in \F_{q}$ since $\F_{q}\cap \F_{q^{r_1}}=\F_{q}$.\\

Equation \ref{array1} implies

\begin{eqnarray*}
\frac{\xi\alpha^{q^{r_1}}+\zeta+\nu\xi}{\alpha^{q^{r_1}}+\nu}=\frac{(a\xi+b)\alpha+a\zeta+a\xi\nu+b\nu}{(c\xi+d)\alpha+c\zeta+c\xi\nu+d\nu}.
\end{eqnarray*}
Recall that $\alpha^{q^{r_1}}=\frac{a\alpha+b}{c\alpha+d}$ so we have, 
\begin{eqnarray*}
\frac{\xi\left(\frac{a\alpha+b}{c\alpha+d}\right)+\zeta+\nu\xi}{\frac{a\alpha+b}{c\alpha+d}+\nu}=\frac{(a\xi+b)\alpha+a\zeta+a\xi\nu+b\nu}{(c\xi+d)\alpha+c\zeta+c\xi\nu+d\nu}.
\end{eqnarray*}
From this, we obtain
\begin{eqnarray*}
\frac{(a\xi+c\xi\nu+c\zeta)\alpha+b\xi+d\nu\xi+d\zeta}{(a+c\nu)\alpha+b+d\nu}=\frac{(a\xi+b)\alpha+a\zeta+a\xi\nu+b\nu}{(c\xi+d)\alpha+c\zeta+c\xi\nu+d\nu}.
\end{eqnarray*}
\
Comparing coefficients, we get
\[
a\xi+c\xi\nu+c\zeta=a\xi+b\] \[b\xi+d\nu\xi+d\zeta=a\zeta+a\xi\nu+b\nu\] \[a+c\nu=c\xi+d\]\[b+d\nu=c\zeta+c\xi\nu+d\nu.
\]
Solving these equations we get \[ \xi=\frac{a+c\nu-d}{c}~\text{and}~\zeta=\frac{b-c\nu\xi}{c}\] where $c\neq 0$.
That is, for any $\nu \in \F_{q}$ we get unique values of $\zeta$ and $\xi$. Since $\nu$ specifies an affine set, this means that each affine set contains exactly one root of $F_{r_1}(x)$. Hence there are $q+1$ roots of $F_{r_1}(x)$ in $O(\alpha)$.
\end{proof}
\end{theorem}
 
Note that if $O(\alpha)$ contains a root of any of the polynomials $F_{r_1}(x)=cx^{q^{r_1}+1}+dx^{q^{r_1}}-ax-b$ it contains precisely $q+1$ roots of the same equation. We know that if $O(\alpha)$ is fixed under $\langle\sigma^{r_1}\rangle$ where $\mbox{gcd}(p,\bar{r}_1)=1$ and $\bar{r}_1\mid (q+1)$ then either $\bar{r}_1=2$ or $\bar{r}_1$ is odd. The case $\bar{r}_1=2$ was addressed in Section \ref{mid}. We have the following lemma.
\begin{lemma}
Suppose $\mbox{gcd}(p,\bar{r}_1)=1$ where $\bar{r}_1$ is an odd integer. Then there are $\frac{\phi(\bar{r}_1)|X(r)|}{2(q+1)}$ projective linear sets fixed by $\langle\sigma^{r_1}\rangle$ if and only if $\bar{r}_1\mid (q+1)$.  
\end{lemma} 

Note that if $\bar{r}_1=r$ then $u=1$ and, by Theorem \ref{main_theo_1}, all the $q+1$ roots of $F_{r_1}(x)$ lie in $\S$. In this case, the number of projective linear sets fixed under $\langle\sigma^{r_1}\rangle$ is $\frac{\phi(\bar{r}_1)}{2}$. We have the following result.
\begin{corollary}\label{Coro_Chika}
Suppose $\mbox{gcd}(p,\bar{r}_1)=1$ and $\bar{r}_1=r$. Then there are $\frac{\phi(\bar{r}_1)}{2}$ projective linear sets fixed under $\langle\sigma^{r_1}\rangle$ if and only if $\bar{r}_1\mid (q+1)$.
\end{corollary} 

Putting all the results together we have proved the following:
\begin{theorem}\label{main_theo_Gamba}
With the notation we have established: 
\begin{enumerate}
\item There are $\frac{|\S|}{q^{3}-q}$ projective linear sets fixed by $\langle\sigma^{r}\rangle$.
\item There are $\frac{\phi(\bar{r}_1)|T(r)|}{2(q-1)}$ projective linear sets fixed by $\langle\sigma^{r_1}\rangle$  if $\bar{r}_1\mid(q-1)$. 
\item If $\mbox{gcd}(\bar{r}_1,r)=p$ then 
\begin{enumerate}
\item there are $\frac{|U(r)|}{q}$ projective linear sets fixed by $\langle\sigma^{r_1}\rangle$ if $\bar{r}_1=p$.
\item there is $1$ projective linear set fixed by $\langle\sigma^{r_1}\rangle$ if $\bar{r}_1=r$.
\end{enumerate}
\item If $\mbox{gcd}(\bar{r}_1,p)=1$ and $\bar{r}_1\mid (q+1)$ then there are 
\begin{enumerate}
\item $\frac{\phi(\bar{r}_1)}{2}$ projective linear sets fixed by $\langle\sigma^{r_1}\rangle$ if $\bar{r}_1=r$.
\item $\frac{\phi(\bar{r}_1)|X(r)|}{2(q+1)}$ projective linear sets fixed by $\langle\sigma^{r_1}\rangle$ if $\bar{r}_1$ is an odd integer.

\end{enumerate}

\end{enumerate}
\end{theorem}
\section{The main theorem}
We can use Theorem \ref{main_theo_Gamba_Aff} and Theorem \ref{main_theo_Gamba} together with the Cauchy-Frobenius theorem to calculate the average number of affine sets and projective  linear sets fixed by an element of $G$ respectively.
%\begin{theorem}
%Let $N(q, r)$ and $N_e(q, r)$ denote the number of orbits in $\S$ under the action of the affine group and the projective linear group and the Frobenius automorphismrespectively, as derived from Theorem \ref{main_theo_Gamba_Aff} and Theorem \ref{main_theo_Gamba}. The number of irreducible and extended irreducible Goppa codes over $\F_q$ of length $q$ and $q+1$ and degree $r$ is at most $N(q, r)$ and $N_e(q, r)$, respectively.
%\end{theorem} 
\begin{theorem}
Let $N(q, r)$ $(resp.~ N_e(q, r))$ denote the number of orbits in $\S$ under the action of the affine group (resp. projective linear group) and the Frobenius automorphism, as derived from Theorem \ref{main_theo_Gamba_Aff} (resp. Theorem \ref{main_theo_Gamba}). The number of irreducible Goppa codes (resp. extended irreducible Goppa codes) over $\F_q$ of length $q$ (resp. $q+1$) and degree $r$ is at most $N(q, r)$ (resp. $N_e(q, r)$).
\end{theorem}

The following table compares $N(q, r)$ and  $N_e(q, r)$.

\begin{table}[htbp]
%\caption{$N(q, r)$ and  $N_e(q, r)$}
%\label{Only_Table}
\centering
\begin{tabular}{cccc}
\hline 
$q$&$r$&$N_e(q, r)$         &$N(q, r)$\\
\hline 
\hline
$2^5$&5&205&6,765 \\ \hline
$5^2$&6&2,667&67,930 \\ \hline
$2^3$&9&29,604&266,304\\ \hline
$3^3$&7&76,027&2,128,684 \\ \hline

\end{tabular}  
\end{table}

\end{document}